\documentclass[12pt,a4paper]{amsart}
\usepackage{amsfonts}
\usepackage{ifthen,verbatim}
\usepackage{amsmath}
\usepackage{graphicx}
\usepackage{amscd,amssymb,amsthm}

\def\be{\begin{equation}}
\def\ee{\end{equation}}

\font\ff=eusm10
\def\K{\hbox{\ff K}}

\newcounter{minutes}
\divide\time by 60
\newcounter{hours}
\multiply\time by 60 \addtocounter{minutes}{-\time}

\title{\bf Landen inequalities for zero-balanced hypergeometric functions}

\author{ Slavko Simi\'c}
\author{Matti Vuorinen}

\address{ Mathematical Institute SANU, Kneza Mihaila 36, 11000
Belgrade, Serbia} \email{ ssimic@turing.mi.sanu.ac.rs}
\address{Department of Mathematics, University of Turku, 20014 Turku,
Finland} \email{vuorinen@utu.fi}

\newtheorem{theorem}[equation]{Theorem}

\newtheorem{lemma}[equation]{Lemma}

\newtheorem{corollary}[equation]{Corollary}
\newtheorem{remark}[equation]{Remark}
\newtheorem{openprob}[equation]{Open problem}

\newtheorem{nonsec}[equation]{}

\numberwithin{equation}{section}

\begin{document}


\begin{abstract}
For zero-balanced Gaussian hypergeometric functions $
F(a,b;a+b;x),$ $a,b>0,$ we determine maximal regions of $ab$ plane
where well-known Landen identities for the complete elliptic
integral of the first kind turn on respective inequalities valid
for each $x\in (0,1)$. Thereby an exhausting answer is given to
the open problem from \cite{avv}.
\end{abstract}

\def\thefootnote{}
\footnotetext{ \texttt{\tiny File:~\jobname .tex,
          printed: \number\year-\number\month-\number\day,
          \thehours.\ifnum\theminutes<10{0}\fi\theminutes }
} \makeatletter\def\thefootnote{\@arabic\c@footnote}\makeatother

\maketitle

{\small \sc Keywords.}{ Log-convexity; Hypergeometric functions;
Inequalities.}

{\small \sc 2010 Mathematics Subject Classification.}{
26D06,33C05}

\section{Introduction}

Among special functions, the hypergeometric function has perhaps
the widest range of applications. For instance, several well-known
classes of mathematical physics are particular or limiting cases of it.
For real numbers $a,b$ and $c$ with $c\neq 0,-1,-2, \cdots $, the
Gaussian hypergeometric function is defined by
\be \label{eq1.1}
F(a,b;c;x) := {}_2 F_1(a,b;c;x)=
\sum_{n=0}^{\infty}\frac{(a,n)(b,n)}{(c,n)}\frac{x^n}{n!}
\ee
for $x\in(-1,1)$, where
$$(a,n):=a(a+1)(a+2) \cdots (a+n-1)$$
for $n=1,2,\cdots$, and $(a,0)=1$ for $a\neq 0$.
For many rational triples $(a,b,c)$
the function \eqref{eq1.1} can be expressed in terms of
elementary functions and long lists
of such particular cases are given in \cite{pbm}.

It is clear that small changes of the parameters $a,b,c$ will have
small influence on the value of $F(a,b;c;x)$. In this paper we shall
study to what extent some well-known properties of the complete
elliptic integral of the first kind
\be \label{eq1.2}
{\K}(x) \equiv \frac{\pi}{ 2} F \Bigl (\frac{1}{ 2},\frac{1}{ 2};1;x^2 \Bigr )
 = \int\limits_0^{\pi/2}(1-x^2\sin^2t)^{-1/2}dt, \quad x \in (0,1),
\ee
can be extended to $F(a,b;a+b;x)$ for $(a,b)$ close to
$(1/2, 1/2)$. Recall that $F(a,b;c;r)$ is called
{\sl zero-balanced} if $c = a+b$. In the zero-balanced case, there
is a logarithmic singularity at $r=1$ and Gauss proved the asymptotic
formula
\be \label{eq1.3}
F(a,b;a+b;r)\sim -\frac{1} {B(a,b)} \log(1-r)
\ee
as $r$ tends to $1$, where
\be \label{eq1.4}
B(z,w)\equiv \frac{{\Gamma(z)\Gamma(w)}}{\Gamma(z+w)},\ Re z>0,
Re w>0
\ee
is the classical beta function. Note that $\Gamma(1/2)= \sqrt{\pi}$ and
$B(\frac{1}{2}, \frac{1}{2})=\pi\,,$ see (\cite[Ch. 6]{as}).

Ramanujan found a much sharper
asymptotic formula 
\be \label{eq1.5}
B(a,b)F(a,b;a+b;r)+\log(1-r)=R(a,b)+O((1-r)\log(1-r))
\ee
as $r$ tends to $1$ (see also \cite{ask}.)
Here and in the sequel,
\be
\label{eq1.6}
\left \{
\begin{array}{rl}
   R(a,b) \equiv & -\Psi(a)-\Psi(b)-2\gamma,
   \quad R({1 / 2}, {1 / 2}) = \log 16 , \\
   \Psi(z) \equiv & \displaystyle{ \frac{d}
   {dz}(\log\Gamma(z))=\frac{\Gamma'(z)}{\Gamma(z)},}\ Re z>0,  \\
\end{array} \right .
\ee
and $\gamma$ is the Euler-Mascheroni constant.  Ramanujan's formula
(\ref{eq1.5}) is a particular case of another well-known formula given
in (\cite[15.3.10]{as}).

\bigskip

We shall use in the sequel the following assertion which is a
mixture of Biernacki-Krzyz and related results on the ratio of
formal power series (\cite{avv},\cite {BOR}).

\bigskip

\begin{lemma}\label{l1}
Suppose that the power series $f(x)=\sum_{n\ge 0}\widehat{f}_n
x^n$ and $g(x)=\sum_{n\ge 0}\widehat{g}_n x^n$ have the radius of
convergence $r>0$ and $\widehat{g}_n>0$ for all $n\in
\{0,1,2,\dots\}$. Denote also
$$
h(x)=\frac{f(x)}{g(x)}=\sum_{n\ge 0}\widehat{h}_n x^n.
$$

1. \ If the sequence $\{\widehat{f}_n/\widehat{g}_n\}_{n\ge 0}$ is
monotone increasing then $h(x)$ is also monotone increasing on
$(0,r)$.

2. \ If the sequence $\{\widehat{f}_n/\widehat{g}_n\}_{n\ge 0}$ is
monotone decreasing then $h(x)$ is also monotone decreasing on
$(0,r)$.

3. \ If the sequence $\{\widehat{f}_n/\widehat{g}_n\}$ is monotone
increasing (decreasing) for $0<n\le n_0$ and monotone decreasing
(increasing) for $n>n_0$, then there exists $x_0\in (0,r)$ such
that $h(x)$ is increasing (decreasing) on $(0,x_0)$ and decreasing
(increasing) on $(x_0,r)$.
\end{lemma}

\bigskip

Some of the most important properties of the
elliptic integral  ${\K}(r)$ are the Landen identities  \cite[p.507]{ww}:
\be \label{eq1.7}
{\K}\Bigl ({ \frac{2\sqrt r}{1+r}}\Bigr)=(1+r){\K(r)},\
{\K}\Bigl ({\frac{1-r}{1+r}}\Bigr)={\frac{1+r}{2}}{\K'}(r),
\ee
where $ \K'(r) = \K(\sqrt{1-r^2}) , r \in (0,1) .$
In \cite[p.79]{avv}, the following problem was raised:

\begin{openprob} \label{prob1.8}
Find an analog of Landen's transformation formulas in
(\ref{eq1.7}) for $F(a,b;a+b;r)$. In particular, if
$k(r)=F(a,b;a+b;r^2)$ and $a,b\in(0,1)$, is it true that
$$k(2\sqrt r/(1+r))\le C k(r)$$
for some constant $C$ and all $r \in (0,1)$?
\end{openprob}

Since $2\sqrt r/(1+r)>r$ for $r\in (0,1)$, $C$ must be greater than 1.

In \cite[pp. 20-21]{avv} and \cite[Theorem 1.4]{abrvv} Gauss' asymptotic formula $(\ref{eq1.3})$
was refined by finding the lower and upper bounds for
$$W(r) = B(a,b)F(a,b;a+b;r)+(1/x)\log(1-r)\,,$$
when $a,b\in(0,1)$ or $a,b\in (1,\infty)$. Our second result
gives a full solution to the Open Problem \ref{prob1.8}.
\bigskip

We wish to point out that in \cite[Thm 1.2(1)]{qv} it was claimed
that for $a,b \in(0,1), c=a+b\le 1,$ the function
\begin{equation}
s(r)= (1+\sqrt{r}) F(a,b;c;r) - F(a,b;c;4
\sqrt{r}/((1+\sqrt{r})^2))
\end{equation}
is increasing in $r \in(0,1)$. As pointed out by A. Baricz \cite{b} the proof
contains a gap and the correct proof will be given here.

We also found another area in $ab$ plane where the function $s(r)$
is monotone decreasing in $r\in(0,1)$.

\section{Main results}
\bigskip

Our first result shows that Landen inequalities hold not only in
the neighborhood of the point $a=b=1/2$ but also in some unbounded
parts of $ab$ plane.

\bigskip

\begin{theorem} \label{s1thm}
For all $a,b>0$ with $ab\le 1/4$ we have that the inequality
\[
F(a,b;a+b;4r/(1+r)^2)\le (1+r)F(a,b;a+b;r^2),
\]
holds for each $r\in(0,1)$.
\bigskip
Also, for $a,b>0, \ 1/a+1/b\le 4$, the reversed inequality
\[
F(a,b;a+b;4r/(1+r)^2)\ge (1+r)F(a,b;a+b;r^2),
\]
takes place for each $r\in (0,1)$.

\bigskip

In the remaining region $a,b>0 \bigwedge ab>1/4 \bigwedge
1/a+1/b>4$ neither of the above inequalities hold for each $r\in
(0,1)$.
\end{theorem}

\bigskip

\begin{figure}
\centering
\includegraphics{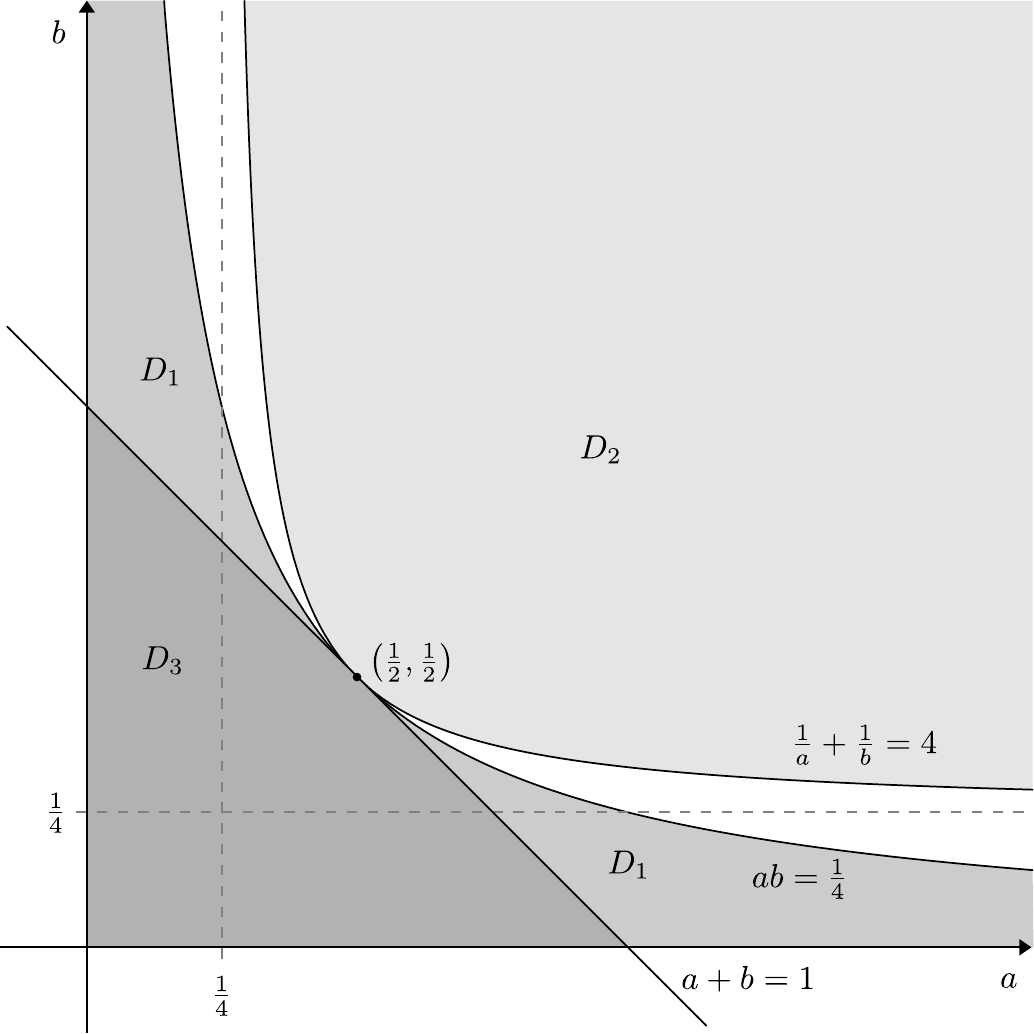}
\caption{The domains $D_j, j=1,2,3$ visualized.   \label{lifig}}
\end{figure}

The disjoint regions in $ab$ plane $D_1=\{(a,b)|a,b>0, \ ab\le
1/4\}$ and \ $D_2=\{(a,b)|a,b>0, 1/a+1/b\le 4\}$, where Landen
inequalities hold, are shown on the Figure \ref{lifig}.

\bigskip

The only common point of the graphs in Figure \ref{lifig} is $(1/2,1/2)$ where equality sign holds.

\bigskip

Two-sided bounds for the ratio of target functions are also
possible.

\bigskip

\begin{theorem}\label{t1}
For each $r\in (0,1)$ and $(a,b)\in D_1$, we have
$$
1<\frac{(1+r)F(a,b;a+b;r^2)}{F(a,b;a+b;4r/(1+r)^2)}<\frac{B(a,b)}{\pi}.
$$

For $(a,b)\in D_2$ the inequalities are reversed,
$$
\frac{B(a,b)}{\pi}<\frac{(1+r)F(a,b;a+b;r^2)}{F(a,b;a+b;4r/(1+r)^2)}<
1.
$$

\end{theorem}

\bigskip

Some numerical estimations of the constant $C$ in Open Problem \ref{prob1.8}
follows.

\begin{corollary}
Let $k(\cdot)$ be defined as in the Open Problem \ref{prob1.8}.
Then, for each $r\in(0,1)$ and $(a,b)\in D_1$, we have
$$
\frac{\pi}{B(a,b)} k(r)<k(2\sqrt{r}/(1+r))<2k(r)
$$

In the region $D_2$ we have
$$
k(r)<k(2\sqrt{r}/(1+r))<\frac{2\pi}{B(a,b)}k(r).
$$
\end{corollary}

\bigskip

Two-sided bounds for the difference exist in a smaller region
$D_3\subset D_1$ (see the picture), where $D_3=\{(a,b)| a,b>0,
a+b\le 1\}$ and in $D_2$.

\bigskip

\begin{theorem}\label{t2}
Let $B=B(a,b)$ be the classical Beta function and $R=R(a,b)$ be
defined by \ref{eq1.6}.

For $a,b>0, a+b\le 1$, we have
$$
0\le (1+\sqrt{r}) F(a,b;a+b;r) - F(a,b;a+b;4
\sqrt{r}/((1+\sqrt{r})^2)) \le (R-\log 16)/B.
$$

If $a,b>0, 1/a+1/b\le 4$, then
$$
0\le  F(a,b;a+b;4 \sqrt{r}/((1+\sqrt{r})^2))-(1+\sqrt{r}) F(a,b;
a+b;r)
 \le (\log 16-R)/B.
$$
\end{theorem}

\bigskip

The second Landen identity has the following counterpart for
hypergeometric functions. The resulting inequalities might be
called Landen inequalities for zero-balanched hypergeometric
functions.

\bigskip

\begin{theorem}\label{t3}
Let $F(x)=F(a,b;a+b;x)$.

For $(a,b)\in D_1$ and each $x\in (0,1)$, we have
$$
\frac{1}{2}<\frac{F((\frac{1-x}{1+x})^2)}{(1+x)F(1-x^2)}<\frac{B(a,b)}{2\pi}.
$$

If $(a,b)\in D_3$, then
$$
(1+x)F(1-x^2)\le 2 F((\frac{1-x}{1+x})^2)\le
(1+x)[F(1-x^2)+(R-\log16)/B].
$$

For $(a,b)\in D_2$, we have
$$
\frac{B(a,b)}{2\pi}
<\frac{F((\frac{1-x}{1+x})^2)}{(1+x)F(1-x^2)}<\frac{1}{2},
$$

and
$$
0\le (1+x)F(1-x^2)-2F((\frac{1-x}{1+x})^2)\le(1+x)(\log16-R)/B.
$$

\end{theorem}

\section{Proofs}
\bigskip

Throughout  this section we denote
$$F(x)=F(a,b;a+b; x), \
G(x)=F(a,b;a+b+1;x),$$
 where $a,b, (a,b)\neq (1/2,1/2)$ are fixed
positive parameters and
$$F_0(x)=F(1/2,1/2;1;x), \
G_0(x)=F(1/2,1/2;2;x),$$
 with the regions $D_1, D_2, D_3$ defined
as above.

\bigskip

 The basic results, which makes possible
all proofs in the sequel, are contained in the following

\begin{lemma}\label{l2} \ 1. \ The function $f(r)=F(r)/F_0(r)$ is monotone decreasing
in $r\in(0,1)$ on $D_1$ and monotone increasing on $D_2$.

2. \ The function $g(r)=G(r)/G_0(r)$ is monotone decreasing on
$D_3$ and monotone increasing on $D_2$.
\end{lemma}

\bigskip

\begin{proof}
We shall use Lemma \ref{l1} in the proof.

\bigskip

Since $\widehat{F}_n=(a)_n(b)_n/(a+b)_n(1)_n, \
\widehat{F_0}_n=((1/2)_n/(1)_n)^2$, applying the lemma one can see
that the monotonicity of $\{\widehat{F}_n/\widehat{F_0}_n\}$
depends on the sign of
\begin{equation}\label{eq1}
T_n=T(a,b;n)=n(ab-1/4)+ab-(a+b)/4=C_1n +C_2.
\end{equation}

Since $(a,b)\neq (1/2,1/2)$ and
$$
C_2=\frac{\sqrt{ab}}{\sqrt{ab}+1/2}C_1-\frac{(\sqrt{a}-\sqrt{b})^2}{4},
$$
it follows

1. \ If $C_1\le 0$ i.e. $(a,b)\in D_1$, then $C_2<0$; hence
$T_n<0$ for $ n=0,1,2,\dots$ and $f(r)$ is monotone decreasing in
$r\in (0,1)$;

2. \ If $C_2\ge 0$ i.e. $(a,b)\in D_2$ then $C_1>0$, that is
$T_n>0, \ n=0,1,2,\cdots$ and $f(r)$ is monotone increasing in
$r$.

\bigskip

In the second case we have
$\widehat{G}_n=(a)_n(b)_n/(a+b+1)_n(1)_n, \
\widehat{G_0}_n=((1/2)_n/(1)_n)^2/(n+1)$ and, proceeding
analogously, we get
$$
T_n=n(ab+a+b-5/4)+2ab-(a+b)/4-1/4=C_3n +C_4.
$$

3. \ If $(a,b)\in D_3$, that is $a,b>0, a+b\le 1$, let $a+b=k>0$.
Then $ab\le k^2/4$ and
$$
C_3\le k^2/4+k-5/4=(k-1)(k+5)/4; \ C_4\le
k^2/2-k/4-1/4=(k-1)(2k+1)/4.
$$

Since $0<k\le 1$, it follows that both $C_3, C_4$ are
non-positive. Therefore $T_n<0, \ n=0,1,2,\dots$ because both
constants cannot be zero simultaneously. By Lemma \ref{l1}, we
conclude that the function $g(r)$ is monotone decreasing in $r\in
(0,1)$.

4. \ If $(a,b)\in D_2$, i.e., $a,b>0, 1/a+1/b\le 4$, then $4ab\ge
a+b\ge 2\sqrt{ab}$, hence $ab \ge 1/4$. Also $a+b\ge 2\sqrt{ab}\ge
2\cdot (1/2)=1$. Therefore $C_3\ge 0$ and
$C_4=(ab-1/4)+(4ab-a-b)/4\ge 0$. As above, we conclude that
$T_n>0, \ n=0,1,2,\cdots$ and $g(r)$ is monotone increasing in
this case.

\end{proof}

\bigskip

\begin{nonsec}{\bf Proof of Theorem \ref{s1thm}.}
{\rm By the above lemma, for each $0<x<y<1$ we have $f(x)>f(y)$ on
$D_1$ and $f(x)<f(y)$ on $D_2$.

Putting $x=x(r)=r^2, \ y=y(r)=4r/(1+r)^2$, we get on $D_1$,
$$
\frac{F(r^2)}{F_0(r^2)}>\frac{F(y)}{F_0(y)},
$$
that is, by Landen's identity,
$$
F(y)<\frac{F_0(y)}{F_0(r^2)} F(r^2)=(1+r)F(r^2).
$$

The second inequality is obtained analogously.

\bigskip

It is easily seen by \eqref{eq1} that in the remaining region the
sequence $\{\widehat{F}_n/\widehat{F_0}_n\}$ decreases and then
increases. By Lemma \ref{l1}, part 3, this means that the function
$f(r)$, for some $r_0\in(0,1)$, decreases in $(0,r_0)$ and
increases in $(r_0,1)$. Therefore, putting $0<x(r)<y(r)<r_0$ and
$r_0<x(r)<y(r)<1$, one concludes that neither of given
inequalities hold for each $r\in (0,1)$.

}\hfill$\square$

\end{nonsec}

\bigskip

\begin{nonsec}{\bf Proof of Theorem \ref{t1}.} \ {\rm Since $f(r)$ is
monotone decreasing on $D_1$, applying Gauss formula, we obtain
$$
1=\lim_{r\to 0^+}\frac{F(r)}{F_0(r)}> \frac{F(r)}{F_0(r)}>
\lim_{r\to
1^-}\frac{F(r)}{F_0(r)}=\frac{B(1/2,1/2)}{B(a,b)}=\frac{\pi}{B(a,b)}.
$$

Therefore,

$$
\frac{F(y(r))}{F(x(r))}<\frac{B(a,b)}{\pi}
\frac{F_0(y(r))}{F_0(x(r))} =(1+r)\frac{B(a,b)}{\pi},
$$
by the Landen identity.

The inequality valid on $D_2$ can be proved
similarly.}\hfill$\square$
\end{nonsec}

\bigskip

\begin{nonsec}{\bf Proof of Theorem \ref{t2}.} \ {\rm Both
assertions of this theorem are a consequence of the following

\begin{lemma}\label{l3} The function
$$
s(r)=(1+\sqrt{r}) F(a,b;a+b;r) - F(a,b;a+b;4
\sqrt{r}/((1+\sqrt{r})^2))
$$
is monotone increasing in $r\in (0,1)$ on $D_3$ and monotone
decreasing on $D_2$.

\end{lemma}}

\begin{proof}
Let $z=\frac{4 \sqrt{r}}{(1+\sqrt{r})^2}$. Then
$$
1-z=\frac{(1-\sqrt{r})^2}{(1+\sqrt{r})^2}; \
\frac{dz}{dr}=\frac{2(1-\sqrt{r})}{\sqrt{r}(1+\sqrt{r})^3}.
$$
Hence
$$
s_1(r):=2\sqrt{r}(1-\sqrt{r})s'(r)=(1-\sqrt{r})F(a,b;a+b;r)+2\sqrt{r}(1-r)F'(a,b;a+b;r)
$$
$$
-\frac{4}{1+\sqrt{r}}(1-z)F'(a,b;a+b;z)
$$
$$
=(1-\sqrt{r})F(a,b;a+b;r)+2\frac{ab}{a+b}\sqrt{r}F(a,b;a+b+1;r)
-\frac{4ab}{(a+b)(1+\sqrt{r})}F(a,b;a+b+1;z)
$$
$$
=(1-\sqrt{r})F(r)+2\frac{ab}{a+b}\sqrt{r}G(r)
-\frac{4ab}{(a+b)(1+\sqrt{r})}G(z).
$$
We used here the well-known formula
\begin{equation}\label{eq3}
(1-x)F'(a,b;a+b;x)=\frac{ab}{a+b}F(a,b;a+b+1;x).
\end{equation}

On the other hand, differentiating the first Landen identity we
get
\begin{equation}\label{eq4}
\frac{1}{1+\sqrt{r}}G_0(z)=
(1-\sqrt{r})F_0(r)+\frac{1}{2}\sqrt{r}G_0(r).
\end{equation}

Since $g(r)$ is monotone decreasing on $D_3$ and $0<r<z<1$, we get
$g(r)>g(z)$ i.e.,
$$
G(z)<\frac{G_0(z)}{G_0(r)} G(r).
$$

This, together with \eqref{eq4}, yields

$$
s_1(r)> (1-\sqrt{r})F(r)+2\frac{ab}{a+b}\sqrt{r}G(r)
-\frac{4ab}{(a+b)(1+\sqrt{r})} \frac{G_0(z)}{G_0(r)} G(r)
$$
$$
= (1-\sqrt{r})F(r)+2\frac{ab}{a+b}\sqrt{r}G(r) -\frac{4ab}{(a+b)}
((1-\sqrt{r})\frac{F_0(r)}{G_0(r)}+\frac{1}{2}\sqrt{r}) G(r)
$$
$$
=(1-\sqrt{r})(F(r)-\frac{4ab}{(a+b)}\frac{F_0(r)}{G_0(r)}G(r)).
$$

By \eqref{eq3} again, we get
$$
\frac{4ab}{(a+b)}\frac{G(r)}{G_0(r)}=\frac{F'(r)}{F_0'(r)}.
$$

Hence,
$$
2\sqrt{r}s'(r)>F(r)- \frac{F'(r)}{F_0'(r)}F_0(r)
=\frac{F^2(r)}{F_0'(r)}\Bigl(\frac{F_0(r)}{F(r)}\Bigr)'.
$$
The last expression is positive on $D_3$ because $D_3\subset D_1$
and, by \eqref{l2}, the function $f(r)= \frac{F(r)}{F_0(r)}$ is
monotone decreasing on $D_1$.

\bigskip

Therefore we proved that the function $s(r)$ is monotone
increasing in $r\in(0,1)$ on $D_3$.

\bigskip

\begin{remark}
Due to the remark in Introduction, this proof gives an affirmative
answer to the 12 years old hypothesis risen in \cite {qv}.

\end{remark}

\bigskip

Since $g(r)$ is increasing on $D_2$, we get
$$
G(z)>\frac{G_0(z)}{G_0(r)} G(r).
$$

Hence, proceeding as before, it follows that
$$
2\sqrt{r}s'(r)<
\frac{F^2(r)}{F_0'(r)}\Bigl(\frac{F_0(r)}{F(r)}\Bigr)'<0,
$$
since $f(r)= \frac{F(r)}{F_0(r)}$ is monotone increasing on $D_2$.

Therefore $s(r)$ is monotone decreasing in $r\in(0,1)$ on $D_2$
and the proof of Lemma \ref{l3} is done.
\end{proof}

\bigskip

{\rm By Lemma \ref{l3} we obtain $\lim_{r\to 0^+}
s(r)<s(r)<\lim_{r\to 1^-} s(r)$ on $D_3$ and $\lim_{r\to 1^-}
s(r)<s(r)<\lim_{r\to 0^+} s(r)$ on $D_2$.

Evidently, $\lim_{r\to 0^+} s(r)=0$.

Applying Ramanujan formula \eqref{eq1.5}, we get
$$
\lim_{r\to 1^-} s(r)=\lim_{r\to 1^-}
(R-2\log(1-r)+\log(1-z)+o(1))/B
$$
$$
=\lim_{r\to
1^-}(R-2\log(1-\sqrt{r})(1+\sqrt{r})+2\log\frac{1-\sqrt{r}}{1+\sqrt{r}}+o(1))/B=(R-\log16)/B.
$$

The assertion of Theorem \ref{t2} follows.} \hfill $\square$

\end{nonsec}

\bigskip

\begin{nonsec}{\bf Proof of Theorem \ref{t3}.} \ {\rm Changing
variable $\frac{1-r}{1+r}=x\in(0,1)$, we obtain
$$
r=\frac{1-x}{1+x}; \ 1+r=\frac{2}{1+x}; \
\frac{4r}{(1+r)^2}=1-x^2.
$$

Putting this in Theorems \ref{t1}, \ \ref{t2}, we obtain
the assertions of Theorem \ref{t3}.} \hfill$\square$

\end{nonsec}

\subsection*{Acknowledgments}
The research of Matti Vuorinen was supported by the Academy of Finland,
Project 2600066611.

\end{document}